\documentclass[journal,10pt]{IEEETran}
\makeatletter
\def\ps@headings{%
\def\@oddhead{\mbox{}\scriptsize\rightmark \hfil \thepage}%
\def\@evenhead{\scriptsize\thepage \hfil \leftmark\mbox{}}%
\def\@oddfoot{}%
\def\@evenfoot{}}
\makeatother
\usepackage[dvips]{graphicx}
\usepackage{subfigure}
\usepackage{amsxtra}
\usepackage{algorithm}
\usepackage{soul}
\usepackage{subfig}
\usepackage{algorithm, algorithmic}
\usepackage{latexsym,url,graphicx,amsmath,amssymb}
\usepackage{epsfig}
\usepackage{tabls}
\usepackage[dvips]{color}
\usepackage{times}
\usepackage{tabularx}
\usepackage{xspace}
\usepackage{pstricks}
\usepackage{multirow}

\usepackage{theorem}

\newtheorem{theorem}{Theorem}[section]
\newtheorem{lemma}{Lemma}[section]

\newtheorem{corollary}{Corollary}[section]

\newenvironment{definition}[1][Definition]{\begin{trivlist}
\item[\hskip \labelsep {\bfseries #1}]}{\end{trivlist}}
\newcommand{\BE}{\begin{equation}}
\newcommand{\EE}{\end{equation}}

\begin{document}

%\title{\LARGE{A Generalized Algorithmic Formulation of Energy and Mutual Information Accumulation in Cooperative Wireless Networks} }

\title{\LARGE{Algorithmic Aspects of Energy-Delay Tradeoff 
in Multihop  \newline Cooperative Wireless Networks} }

%{On Efficient Energy and Information Distribution in Cooperative Multihop Wireless Networks}

%On Energy Efficient Transmission in Cooperative Multihop Wireless Networks}

%Joint Power Control and Scheduling in Cooperative Wireless Networks}

\author{\IEEEauthorblockN{Marjan Baghaie, \emph{\small Student Member, IEEE}, Bhaskar Krishnamachari, \emph{\small Member, IEEE}, Andreas F. Molisch, \emph{\small Fellow, IEEE}}%\IEEEauthorblockN{Marjan Baghaie, and Bhaskar Krishnamachari}\\
%\IEEEauthorblockA{Department of Electrical Engineering, Viterbi School of Engineering\\
%University of Southern California, Los Angeles, California 90089\\
%Email: \{baghaiea, bkrishna\}@usc.edu}}
\thanks{The authors are with the University of Southern California, Los Angeles, CA, 90007. Email: \{baghaiea, bkrishna, molisch\}@usc.edu.

This research was sponsored in part by the U.S. Army Research
Laboratory under the Network Science Collaborative Technology
Alliance, Agreement Number W911NF-09-2-0053, and by U.S. National
Science Foundation under awards CNS-0627028 and CNS-1049541.}
}

% make the title area
\maketitle
\begin{abstract}
We consider the problem of energy-efficient transmission in cooperative multihop wireless networks. Although the performance gains of cooperative approaches are well known, the combinatorial nature of these schemes makes it difficult to design efficient polynomial-time algorithms for deciding which nodes should take part in cooperation, and when and with what power they should transmit. In this work, we tackle this problem in \emph{memoryless} networks with or without delay constraints, i.e., quality of service guarantee. We analyze a wide class of setups, including unicast, multi-cast, and broadcast, and two main cooperative approaches, namely:
energy accumulation (EA) and mutual information accumulation (MIA).
We provide a generalized algorithmic formulation of the problem that encompasses all those cases. We investigate the similarities and differences of EA and MIA in our generalized formulation.  We prove that the \emph{broadcast and multicast} problems are, in general, not only NP hard but also $o(\log(n))$ inapproximable. We break these problems into three parts: ordering, scheduling and power control, and propose a novel algorithm that, given an ordering, can optimally solve the joint power allocation and scheduling problems
simultaneously in polynomial time. We further show empirically that this algorithm used in conjunction with an ordering derived heuristically using the
Dijkstra's shortest path algorithm yields near-optimal performance in typical settings. For the \emph{unicast} case, we prove that although the problem remains NP complete with MIA, it can be solved optimally and in polynomial time when EA is used. We further use our algorithm to study numerically the
trade-off between delay and power-efficiency in cooperative broadcast and compare the performance of EA vs MIA as well as the performance of our cooperative algorithm with a smart non-cooperative algorithm in a broadcast setting.  
\end{abstract}
%\begin{keywords}
%resource allocation, scheduling, fountain codes, Cooperative transmission, wireless networks, broadcast, multicast, optimization
%time.
%\end{keywords}
\IEEEpeerreviewmaketitle

\section{INTRODUCTION}

In a wireless network, a transmit signal intended for one node is received not only by that node but also by other nodes. In a traditional point-to-point system, where there is only one intended recipient, this innate property of the wireless propagation channel can be a drawback, as the signal constitutes undesired interference in all nodes but the intended recipient. However, this effect also implies that a packet \emph{can} be transmitted to multiple nodes simultaneously without additional energy expenditure. Exploiting this Òbroadcast advantageÓ, broadcast, multicast and multihop unicast systems can be designed to work cooperatively and thereby achieve potential performance gains. As such, cooperative transmission in wireless networks has attracted a lot of interest not only from the research community in recent years~\cite{Khandani03, Maric04, Mergen06, Chen05, Mergen07, Jakllari07, Draper08} but also from industry in the form of first practical cooperative mobile ad-hoc network systems~\cite{Halford10}.

We focus on the problem of cooperative transmission in this work, where a single node is sending a packet to either the entire network (broadcast), a single destination node (unicast) or more than one destination node (multicast), in a multihop wireless network. Other nodes in the network, that are neither the source nor the destination, may act as relays to help pass on the message through multiple hops. The transmission is completed when all the destination nodes have successfully received the message. We particularly focus on the case where there is a delay constraint, whereby the destination node(s) should receive the message within the delay constraint, however, we also discuss how our results apply to the unconstrained case. 

A key problem in such cooperative networks is routing and resource allocation, i.e., the question which nodes should participate in the transmission of data, and when, and with how much power, they should be transmitting. The situation is further complicated by the fact that the routing and resource allocation depends on the type of cooperation and other details of the transmission/reception strategies of the nodes. We consider in this paper a time-slotted system in which the nodes that have received and decoded the packet are allowed to re-transmit it in future slots. During reception, nodes add up the signal power (energy accumulation, EA) or the mutual information (mutual information accumulation, MIA) received from multiple sources. EA, which has been discussed in prior work~\cite{Maric04, Chen05, Mergen06, Mergen07}, can be implemented by using maximal ratio combining (MRC) of orthogonal signals from source nodes that use orthogonal time/frequency channels, or spreading codes, or distributed space-time codes. MIA can be achieved using rateless codes~\cite{Draper08, Molisch07}. Although these techniques are often treated separately in the literature, we shall see how our formulation of the problem encompasses both approaches and allows many of the results to be extended to both. 

We furthermore assume that the nodes are memoryless, i.e., accumulation at the receiver is restricted to transmissions from multiple nodes in the present time slot, while signals from previous timeslots are discarded. This assumption is justified by the limited storage capability of nodes in ad-hoc networks, as well as the additional energy consumption nodes have to expand in order to stay in an active reception mode when they ÒoverhearÓ weak signals in preceding timeslots. Note that much of the literature cited above has used the assumption of nodes with memory, so that their results are not directly comparable to ours.

A key tradeoff is between the total energy consumption\footnote{As we consider fixed time slot durations, we use the words energy and power interchangeably in this paper} and the total delay measured in terms of the number of slots needed for all destination nodes in the network to receive the message. At one extreme, if we wish to minimize delay, each transmitting node should transmit at the highest power possible so that the maximum number of receivers can decode the message at each step (indeed, if there is no power constraint, then the source node could transmit at a sufficiently high power to reach all destination nodes in the first slot itself). On the other hand, reducing transmit power levels to save energy, may result in fewer nodes decoding the signal at each step, and therefore in a longer time to complete the transmission. We therefore formulate the problem of performing this transmission in such a way that the total transmission energy over all transmitting nodes is minimized, while meeting a desired delay constraint on the maximum number of slots that may be used to complete the transmission. The design variable in this problem is to decide which nodes should transmit, when, and with what power.

 \begin{table*} [!htb]
\centering

\begin{tabular}{|c|c|c|c|c|} 
  \hline
 \multirow{2}{*}{ \bf NEGATIVE RESULTS} & \multicolumn{2}{|c|} {\bf Energy Accumulation}   &\multicolumn{2}{|c|} { \bf Mutual Information Accumulation} \\ \cline{2-5}
  & \it {Delay Constraint ($T$)}& \it{ Unconstrained} & \it{Delay Constraint ($T$)} & \it{Unconstrained} \\ \cline{2-5}
    \hline \hline
{ \bf Broadcast} & $o(\log(n))$ inapproximable  for $T\geq 3$ & NP-complete&$o(\log(n))$ inapproximable for $T\geq 3$ & NP-complete\\
  \hline \hline
  { \bf Multicas}t & $o(\log(n))$  inapproximable  for $T\geq 3$& NP-complete&$o(\log(n))$  inapproximable for $T\geq 3$ & NP-complete\\
  \hline \hline
  { \bf Unicast}& Polynomial time & Polynomial time & NP-complete for $T\geq 4$ & --- \\
  \hline
\end{tabular}
\caption{Summary of the algorithmic negative results.}
\label{tb:neg}
\end{table*}

\begin{table*} [!htb]
\centering

\begin{tabular}{|p{3cm}|p{3cm}|p{3cm}|p{3cm}|p{3cm}|} 
  \hline
 \multirow{2}{*}{ \bf POSITIVE RESULTS} & \multicolumn{2}{|c|} {\bf Energy Accumulation}   &\multicolumn{2}{|c|} { \bf Mutual Information Accumulation} \\ \cline{2-5}
  & \it {Delay Constraint ($T$)}& \it{ Unconstrained} & \it{Delay Constraint ($T$)} & \it{Unconstrained} \\ \cline{2-5}
    \hline \hline
{ \bf Broadcast} & $\bullet O(n^\epsilon)$ for $\epsilon > 0$ \newline $\bullet O(T \log^2(n))$, for fixed $T$ \newline $\bullet$ Polynomial time given ordering (DMCT\_go)  \newline $\bullet$ Polynomial time algorithm for $T=1,2$ &Polynomial time given ordering \newline $1$D dynamic program given in $(7)$&  $\bullet$ Polynomial time given ordering (DMCT\_go) \newline $\bullet$ Polynomial time algorithm for $T=1,2$ & Polynomial time given ordering \newline $1$D dynamic program given in $(7)$ \\
  \hline \hline
  { \bf Multicas}t & Polynomial time given ordering (DMCT\_go) & Polynomial time given ordering (DMCT\_go) & Polynomial time given ordering (DMCT\_go) & Polynomial time given ordering (DMCT\_go) \\
  \hline \hline
  { \bf Unicast}& Polynomial time & Polynomial time & $\bullet$ Polynomial time given ordering (DMCT\_go)  $T\geq 4$ \newline $\bullet$ Polynomial time algorithm for $T=1,2$ & --- \\
  \hline
\end{tabular}
\caption{Summary of the algorithmic positive results.}
\label{tb:pos}
\end{table*}

The key contributions of our work are as follows:
\begin{itemize}

\item We formulate the problem of minimum energy transmission in cooperative networks. Although the prior literature have focused on either EA (\cite{Maric04, Mergen07}) or MIA (\cite{Draper08, Molisch07}) and have treated them separately, our generalized formulation can treat both methods as variations of the same problem.

\item Our formulation of delay-constrained minimum energy broadcast in cooperative networks, goes beyond the prior work in the literature on cooperative broadcast which has focused either on minimizing energy without delay constraints~\cite{Maric04, Mergen07}, or on delay analysis without energy minimization~\cite{Baghaie09}. Our extended problem formulation allows us to expose and investigate the energy-delay tradeoffs inherent in cooperative networking.

\item We not only prove that the delay constrained minimum energy cooperative broadcast (DMECB) and multicast (DMECM) problems are NP-complete, but also that they are $o(\log (n))$-inapproximable (i.e., unless $P=NP$, it is not possible to develop a polynomial time algorithm for this problem that can obtain a solution that is strictly better than a logarithmic-factor of the optimum in all cases). We are not aware of prior work on cooperative broadcast or multicast that shows such inapproximability results.

\item We show that the delay constrained minimum energy cooperative unicast (DMECU) problem is solvable in polynomial time using EA but is NP-complete using MIA. We are unaware of any hardness results on unicast approaches using mutual information accumulation.

\item For the cases where we prove the transmission problem to be hard, we are able to show that for any given \emph{ordering} of the transmissions (which dictates that a node later in the ordering may not transmit before the nodes earlier in the ordering have decoded successfully), then the problem of joint scheduling and power allocation can in fact be solved optimally in polynomial time using a combination of Dynamic Programming for the scheduling and Convex Programming for the power allocation.

\item For small network instances, we compute the optimal solution through exhaustive search, and show empirically through simulations that our proposed joint scheduling and power control method works near-optimally in typical cases when used in conjunction with an ordering provided by the Dijkstra tree construction.

\item We also show through simulations the delay-energy tradeoffs and minimum energy performance for larger networks and demonstrate the significant improvements that can be achieved by our solution compared to non-cooperative broadcast. We further compare the performance of our proposed broadcast algorithm under MIA and EA approaches.

\item For DMECB where EA is used, we present a reduction that would allow for a polynomial time algorithm for the joint ordering-scheduling-power control problem that is provably guaranteed to offer a $O(n{^\epsilon})$ approximation, for any $\epsilon > 0$. This algorithm is based on the current best-known algorithm for the bounded diameter directed Steiner tree problem~\cite{Moses99}. Using the same reduction, we can also get an approximation factor of  $O(T \log^2(n))$ for a fixed delay constraint $T$. Given that DMECB is $o(\log(n))$ inapproximable for any $T > 2$, this provides a fairly tight approximation, especially for small $T$.

\end{itemize}

This paper makes several contributions that significantly enhance our understanding of complexity and algorithm design for cooperative transmission in wireless networks, in the context of energy-delay tradeoffs. The summary of the algorithmic results developed in this paper are presented in Tables \ref{tb:neg} and \ref{tb:pos}. 
It is worth noticing that although in this work we have considered a memoryless setting, all the negative results presented in Table \ref{tb:neg} extend to the case where there is no memory.

The rest of the paper is organized as follows: section~\ref{sec:related} places our work in the context of prior related work. Section ~\ref{sec:model} describes the system model. The problem formulation is presented in section \ref{sec:formul}.  In section \ref{sec:npcomp}, we establish hardness results. The main positive result for the hard problems described in \ref{sec:npcomp} is presented in section~\ref{sec:optBroadcast}, where a polynomial time algorithm for optimum delay constrained scheduling and power allocation is presented, for the case when the ordering is given. In section \ref{sec:unicast}, we provide a polynomial time optimal algorithm for the unicast case where EA is used. Section~\ref{sec:approx} provides approximation results for DMECB when EA is used. Simulation results are presented in section~\ref{sec:simul}, where we suggest and evaluate several heuristics for the ordering. We present concluding comments and directions for future work in section~\ref{sec:conc}, along with two tables summarizing the results.

\section{RELATED WORK \label{sec:related}}

EA and MIA are two of the main approaches\footnote{Other cooperation schemes include distributed beamforming \cite{Khandani03} and coded cooperation \cite{Nosratinia04}, which will not be further considered in this paper.} for cooperative communications in wireless network .

In EA~\cite{Maric04, Chen05, Mergen06, Mergen07}, a receiver can recover the original packet so long as the total received energy from multiple sources or successive transmissions exceeds a given threshold. Such an approach can be implemented using maximal ratio combining of orthogonal signals from multiple sources, e.g., through a Rake receiver in CDMA or distributed space-time codes. It has been shown that one can achieve significant saving in energy and/or transmission time when using an energy accumulation protocol, compared to traditional protocols~\cite{Maric04, Mergen07, Baghaie09}. If energy accumulation is achieved by transmitting the exact same packet either from different relays or through successive re-transmissions, the scheme is shown to achieve capacity in an asymptotically wideband regime~\cite{Maric04}.  We note that recently a commercially developed cooperative mobile ad hoc network system has been developed which utilizes a pragmatic cooperation method requiring minimal information exchange, based on a combination of phase dithering and turbo codes~\cite{Lee06, Halford10}. It is shown in~\cite{Lee06} that the performance of this pragmatic scheme is close to that of an ideal energy-accumulation approach based on space-time coding.

In MIA ~\cite{Draper08, Molisch07}, the receiving node accumulates mutual information for a packet from multiple transmissions until it can be decoded successfully. In practice, this can be achieved using rateless codes~\cite{Draper08, Molisch07}. The two schemes have been shown to be equivalent at low signal-to-noise ratios (SNRs), while in high SNR regimes the latter is shown to have superior performance~\cite{Draper08, Molisch07}.

Many network protocols in mobile ad hoc and sensor networks need to operate in broadcast mode to disseminate certain control messages to the entire network (for instance, to initiate route requests, or to propagate a query). The subject of broadcast transmission in multi-hop wireless networks has attracted a lot of attention from the research community in both non-cooperative~\cite{Ni99, Williams02, Cagalj02} and cooperative settings~\cite{Maric04, Mergen06, Mergen07, Baghaie09, Jakllari07}. For traditional non-cooperative wireless networks, Cagalj \emph{et al.} \cite{ Cagalj02} showed that the problem of minimum energy broadcast is NP-hard. In~\cite{Mergen06}, Mergen et al. show through a continuum analysis the existence of a phase transition in the behavior of cooperative broadcast: if the decoding threshold is below a critical value then the broadcast is successful, else only a fraction of the network is reached. In~\cite{Mergen07}, Mergen and Scaglione, show that the problem of scheduling and power control for minimum energy broadcast is tractably solvable for highly dense (continuum) networks and show the gains obtained with respect to noncooperative broadcast. In~\cite{Baghaie09}, we examined the delay performance of cooperative broadcast and show that cooperation can result in extremely fast message propagation, scaling logarithmically with respect to the network diameter, unlike the linear scaling for non-cooperative broadcast.

The work by Maric and Yates~\cite{Maric04} is closest in spirit to our work. They too address the minimum-energy cooperative broadcast problem where nodes use EA. However, in their work the system has memory, in that the nodes can save soft information from all previous transmissions throughout time and use it to decode data later on.  They prove that the problem is NP-complete in this case. In their setting, because of the memory, it suffices to have each transmitter transmit only once; therefore there is no distinction between ordering and scheduling. This is no longer true in our memoryless setting where the energy from past transmissions cannot be accumulated. A further key distinction in our work is that we consider delay constraints, whereas~\cite{Maric04} focuses only on the minimum energy cooperative solution without delay constraints; additionally we consider MIA, which is not addressed in ~\cite{Maric04}.

One prior work that discusses the power-delay tradeoff in a cooperative setting is \cite{Gold05}; however, the focus of that work is on space-time codes used for unicast, not broadcast. In earlier work~\cite{Baghaie11}, we had considered this tradeoff in a broadcast cooperative setting using EA and conjectured that many of the results would extend to MIA but the investigation of that conjecture had remained an open problem. In this paper, we are addressing that open problem, as well as considering unicast and multicast settings - which were not addressed in ~\cite{Baghaie11}.

MIA is most notably discussed in~\cite{Draper08, Molisch07}. However, the discussions are focused on unicast routing. The algorithms presented are heuristics, the performance of which are verified via simulations. We are not aware of any hardness results on cooperative MIA or any discussions on the use of MIA in a general cooperative broadcast setting.

\section{SYSTEM MODEL \label{sec:model}}

We consider a wireless network with $n$ nodes. Radio propagation is modeled by a given symmetric $n$ by $n$ static channel matrix, $H = \{h_{ij}\}$, representing the (power) gain on the channel between each pair of nodes $i$ and $j$. Time is assumed to be discretized into fixed-duration slots; without loss of generality we assume unit slot durations. We assume cooperative communication in the receivers, encompassing two scenarios: EA and MIA. Only a single message is transmitted through the network. 

In EA, the received power at a given receiver in a specific timeslot is sum of the powers received from the transmitters that are active during that slot. As described in~\cite{Mergen06, Mergen07, Chen05}, this kind of additive received power can be achieved via \emph{maximal ratio combining} under different scenarios including transmission using TDMA, FDMA channels, as well as with CDMA spreading codes and space-time codes. 

MIA can be implemented using rateless codes and decoders at receivers, as described in~\cite{Draper08,Molisch07,Castura07}.  With proper design (e.g., different spreading codes), information streams from different relay nodes can be distinguished, and the mutual information of signals transmitted by different nodes can be accumulated. We consider a per-node bandwidth constraint and dynamic power allocation.

We assume appropriate coding is used so that each receiving node can decode the message so long as its accumulated received mutual information exceeds a given threshold $\theta$ that represents the bandwidth-normalized entropy of the information codeword in nats/Hz. Furthermore, all nodes are assumed to operate in half-duplex mode, i.e. they cannot transmit and receive simultaneously. If used in transmission, the nodes operate based on a decode and forward protocol. Therefore, they are not allowed to take part in transmission until they have fully decoded their message.

Assuming the noise power is the same at all receivers, we can assume without loss of generality the noise power to be normalized to unity so that the transmit power attenuated by the channel becomes equivalent to the signal to noise ratio (SNR). As mentioned in Sec. I, we assume a memory-less model in which nodes do not accumulate energy or information from transmissions occurred in previous time slots. 

\section{PROBLEM FORMULATION}\label{sec:formul}

In this section we provide a generalized formulation for the delay constrained minimum-energy cooperative transmission (DMECT) problem in the setting described in Sec. III. 

We assume that the transmission begins from a single source node. The aim is to get the message to all the nodes in a destination set $\cal{D}$, with the minimum possible total energy, within a time $T$ (which can take on any value from $1$ to $n-1$).  Every node in the network is allowed to cooperate in the transmission, so long as they have already decoded the message. The problem now becomes: which nodes should take part in the cooperation, when and with what power should they transmit to achieve this aim while meeting the constraints and incurring minimum total transmission power.

Recalling the memoryless assumption, the condition for successful decoding at some receiver node $r$ at time $t$ when a set of nodes $S(t)$ is transmitting packets,  with transmit power $p_{st}, \forall s\in S(t)$ is:

\begin{equation} y_{rt}  \geq \theta ~\label{eqn:decode}\end{equation}
with $y_{rt}$ being the mutual information accumulated by node $r$ at time $t$. Let $x_{it}$ be an indicator binary variable that indicates whether or not node $i$ is allowed to transmit at time $t$. In other words, we define $x_{it}$ to be $1$, if node $i$ is allowed to transmit at time $t$ (i.e. has decoded the message by the beginning of time slot $t$ as per equation~(\ref{eqn:decode})), and $0$ otherwise. Let $p_{it}$ be the transmit power for each node $i$ at each time $t$. Without loss of generality, the source node is assigned node index $1$.

The DMECT problem can then be formalized as a combinatorial optimization problem:

\begin{eqnarray} \label{eqn:DMECT_opt} \min & P_{total} = \sum_{t=1}^{T}\sum_{i=1}^n  p_{it} \\  \begin{tabular}{r}
  s.t. \\
 \\
    \\
\\
    \\
\\ \end{tabular}
& \begin{tabular}{ >{$}l<{$} l }
  1. &$ p_{it} \geq 0,~~\forall i , \forall t$   \\
  2. &$ x_{iT+1} \geq 1,~~\forall i \in \cal{D} $  \\
  3. & $x_{it+1} \leq {1\over \theta}y_{it}+x_{it},~~\forall i, \forall t$  \\
  4. &$ x_{1t}=1, ~~\forall t$  \\
   5. &$ x_{i1}=0,~~ \forall i \neq 1 $ \\
6. & $x_{it} \in \{0,1\}  $\\
\end{tabular} \nonumber \end{eqnarray}

where, for the energy accumulation (EA) case\footnote{Notice that because of the monotonicity of the $\log$ function, $y_{it}  \geq \theta$ in this case is equivalent to $\sum_{s \in S(t)} p_{st} x_{st} h_{si} \geq e^\theta -1$}: 

\begin{equation}y_{it}= \log\left(1+\sum_{s \in S(t)} p_{st} x_{st} h_{si}\right) \end{equation}

and for mutual information accumulation (MIA) case:

\begin{equation}y_{it}= \sum_{s \in S(t)} \log \left(1+p_{st} x_{st} h_{si} \right) \end{equation}

Constraint $2$ ensures that every nodes in the destination set successfully decodes the message within the time constraint $T$, constraint $3$ ensures that a node cannot transmit unless it has already received the message while simultaneously making sure that a node that has decoded the message in previous time slots will not be prevented from transmitting in future time slots (if it wants to transmit), constraint $4$ assigns the source node, and all other constraints are self-explanatory.

In general, there are three variations of this problem, based on the size of the destination set:
\begin{itemize}

\item Delay constrained minimum energy cooperative unicast (DMECU): where the set $\cal{D}$ includes a single destination node. 

\item Delay constrained minimum energy cooperative multicast (DMECM): where the set $\cal{D}$ includes more than one destination node. 

\item Delay constrained minimum energy cooperative broadcast (DMECB):  where the set $\cal{D}$ includes all the nodes apart from the source node.

\end{itemize}

The decision version of these problem, can be defined correspondingly as follows: ``Given some power bound $C$, does there exist an allocation of powers, $p_{it}$, satisfying the constraints in (\ref{eqn:DMECT_opt}) such that $P_{total} \leq C$?"  An instance of this decision problem is defined by giving the symmetric $n\times n$ matrix $H$,  with a designated source node (vertex), a  destination set $\cal{D}$, a delay bound $T$, and a power bound $C$.

Notice that assigning $T \geq n$, in the above formulation, results in the problem definition in the case where there is no delay constraint.  Note also that a requirement for per-node maximum power can be trivially added to the above formulation as additional constraint; we have left that out for simplicity. Should the maximum power be added, it should be large enough to ensure a feasible solution exists for the given connectivity and delay constraint.

\section{HARDNESS \label{sec:npcomp}}

In this section we prove that finding an optimal solution for DMECB and DMECM problems is not only NP-hard but also $o(\log(n))$ inapproximable Ð i.e., finding any polynomial time algorithm that approximates the optimal solutions within a factor of $o(\log(n))$ is also NP-hard. We show this by demonstrating that any instance of the set cover problem can be reduced to an instance of DMECB (and by extension DMECM). We further prove NP-completeness for DMECU when MIA is used; note that DMECU with EA will be treated in Sec. \ref{sec:unicast}.

\subsection{Set Cover Problem}
The set cover problem is a classical problem in computer science \cite{Vazi01}. It is stated as follows: Given a universe $U$ of $n$ elements and a collection of subsets of $U$, $S={S_1,S_2,...S_k}$,  find a minimum subcollection of $S$ that covers all elements of $U$. This problem is NP-complete and was shown, in \cite{Raz97}, to be $o(\log(n))$ inapproximable. 

The set cover problem can be thought of as a bipartite graph $G(V,E)$, with $|V|=k+n$, representing the $k$ sets and $n$ elements and the edges are used to connect each set to its elements. This is shown in Figure \ref{fig:setcover} (a), where we assign a vertex for each set in the top part of the graph, and assign a vertex for each element in the bottom part of the graph. We connect each set to its elements using an edge. One can think of each vertex in this graph as a node in a network, in which edges exist between any pair of nodes for which $h_{ij} > 0$, and the edges are labeled with a weight $w_{ij} $ that corresponds to the transmit power needed at node $i$ to exceed a threshold of  $\theta$ at the receiver $j$, in a single time slot if $i$ was the only transmitter. Given an instance, $G$, of the set cover problem, the optimal solution to the set cover problem, $OPT_{sc}$, would find the minimum subset of vertices in the top part of the graph, so that their transmission of a message can broadcast the message to all the vertices in the bottom part of the graph.

\subsection{Inapproximablity of DMECB} \label{sec:npcomB}

\begin{figure}[!h]
\centering

\scalebox{1} % Change this value to rescale the drawing.
{
\begin{pspicture}(0,-2.250039)(8.49,2.290039)
\rput{-90.0}(1.0949609,0.27503905){\pscircle[linewidth=0.04,dimen=outer](0.685,-0.40996096){0.1}}
\rput{-90.0}(1.7149609,0.895039){\pscircle[linewidth=0.04,dimen=outer](1.305,-0.40996096){0.1}}
\rput{-90.0}(2.9349608,2.155039){\pscircle[linewidth=0.04,dimen=outer](2.545,-0.38996094){0.1}}
\rput{-90.0}(2.0749607,-1.6249608){\pscircle[linewidth=0.04,dimen=outer](0.225,-1.8499608){0.1}}
\rput{-90.0}(2.6749609,-1.0249608){\pscircle[linewidth=0.04,dimen=outer](0.825,-1.8499608){0.1}}
\rput{-90.0}(4.954961,1.1750393){\pscircle[linewidth=0.04,dimen=outer](3.065,-1.8899608){0.1}}
\rput{-90.0}(3.3549607,-0.34496075){\pscircle[linewidth=0.04,dimen=outer](1.505,-1.8499608){0.1}}
\psline[linewidth=0.04cm](0.645,-0.46996096)(0.225,-1.7699609)
\psline[linewidth=0.04cm](0.705,-0.46996096)(1.465,-1.7499609)
\psline[linewidth=0.04cm](1.285,-0.46996096)(0.845,-1.7499609)
\psline[linewidth=0.04cm](2.5849218,-0.43003905)(3.005,-1.8299608)
\psline[linewidth=0.04cm](2.485,-0.42996094)(0.905,-1.8099608)
\psdots[dotsize=0.06,dotangle=-270.0](2.065,-1.889961)
\psdots[dotsize=0.06,dotangle=-270.0](2.265,-1.889961)
\psdots[dotsize=0.06,dotangle=-270.0](2.465,-1.889961)
\psframe[linewidth=0.04,dimen=outer](3.244922,-0.09003905)(0.0449219,-2.250039)
\rput{-90.0}(3.514961,7.095039){\pscircle[linewidth=0.04,dimen=outer](5.305,1.7900391){0.1}}
\rput{-90.0}(3.594961,5.2150393){\pscircle[linewidth=0.04,dimen=outer](4.405,0.81003904){0.1}}
\rput{-90.0}(4.214961,5.835039){\pscircle[linewidth=0.04,dimen=outer](5.025,0.81003904){0.1}}
\rput{-90.0}(5.434961,7.095039){\pscircle[linewidth=0.04,dimen=outer](6.265,0.830039){0.1}}
\psline[linewidth=0.04cm](5.245,1.7900391)(4.445,0.91003907)
\psline[linewidth=0.04cm](5.285,1.730039)(5.045,0.890039)
\psline[linewidth=0.04cm](5.365,1.730039)(6.205,0.91003907)
\usefont{T1}{ptm}{m}{n}
\rput(5.31,2.0950391){\small $r$}
\psline[linewidth=0.04cm,linestyle=dashed,dash=0.16cm 0.16cm](3.725,1.1100391)(7.525,1.1100391)
\usefont{T1}{ptm}{m}{n}
\rput(7.73,1.725039){\small Level $0$}
\usefont{T1}{ptm}{m}{n}
\rput(7.73,0.88066405){\small Level $1$}
\usefont{T1}{ptm}{m}{n}
\rput(7.75,-0.27386716){\small Level $2$}
\usefont{T1}{ptm}{m}{n}
\rput(7.71,-1.8138671){\small Level $3$}
\rput{-90.0}(4.814961,3.995039){\pscircle[linewidth=0.04,dimen=outer](4.405,-0.40996096){0.1}}
\rput{-90.0}(5.434961,4.615039){\pscircle[linewidth=0.04,dimen=outer](5.025,-0.40996096){0.1}}
\rput{-90.0}(6.654961,5.875039){\pscircle[linewidth=0.04,dimen=outer](6.265,-0.38996094){0.1}}
\psline[linewidth=0.04cm](4.384922,0.70996094)(4.365,-0.32996106)
\psline[linewidth=0.04cm](5.044922,0.729961)(5.025,-0.30996084)
\psline[linewidth=0.04cm](6.2449217,0.74996096)(6.2449217,-0.31003904)
\rput{-90.0}(5.794961,2.0950394){\pscircle[linewidth=0.04,dimen=outer](3.945,-1.8499608){0.1}}
\rput{-90.0}(6.394961,2.6950393){\pscircle[linewidth=0.04,dimen=outer](4.545,-1.8499608){0.1}}
\rput{-90.0}(8.674961,4.895039){\pscircle[linewidth=0.04,dimen=outer](6.785,-1.8899608){0.1}}
\rput{-90.0}(7.0749607,3.3750393){\pscircle[linewidth=0.04,dimen=outer](5.225,-1.8499608){0.1}}
\psline[linewidth=0.04cm](4.365,-0.46996096)(3.945,-1.7699609)
\psline[linewidth=0.04cm](4.425,-0.46996096)(5.185,-1.7499609)
\psline[linewidth=0.04cm](5.005,-0.46996096)(4.565,-1.7499609)
\psline[linewidth=0.04cm](6.304922,-0.47003904)(6.725,-1.8299608)
\psline[linewidth=0.04cm](6.205,-0.42996094)(4.625,-1.8099608)
\psdots[dotsize=0.06,dotangle=-270.0](5.365,0.44003916)
\psdots[dotsize=0.06,dotangle=-270.0](5.565,0.44003916)
\psdots[dotsize=0.06,dotangle=-270.0](5.765,0.44003916)
\usefont{T1}{ptm}{m}{n}
\rput(4.15,0.33785155){\small $M$}
\psline[linewidth=0.04cm,linestyle=dashed,dash=0.16cm 0.16cm](3.685,0.11003905)(7.485,0.11003905)
\psline[linewidth=0.04cm,linestyle=dashed,dash=0.16cm 0.16cm](3.725,-1.5299609)(7.525,-1.5299609)
\psframe[linewidth=0.04,dimen=outer](6.944922,-0.09003905)(3.764922,-2.250039)
\psdots[dotsize=0.06,dotangle=-270.0](5.785,-1.8999609)
\psdots[dotsize=0.06,dotangle=-270.0](5.985,-1.8999609)
\psdots[dotsize=0.06,dotangle=-270.0](6.185,-1.8999609)
\psline[linewidth=0.04cm](3.784922,-1.070039)(6.944922,-1.070039)
\psdots[dotsize=0.06,dotangle=-270.0](5.425,-0.42996106)
\psdots[dotsize=0.06,dotangle=-270.0](5.625,-0.42996106)
\psdots[dotsize=0.06,dotangle=-270.0](5.825,-0.42996106)
\usefont{T1}{ptm}{m}{n}
\rput(4.83,0.33785155){\small $M$}
\usefont{T1}{ptm}{m}{n}
\rput(6.45,0.35785156){\small $M$}
\usefont{T1}{ptm}{m}{n}
\rput(3.9,1.475039){\small $G'$}
\psline[linewidth=0.04cm](3.224922,-1.0300391)(0.0449219,-1.070039)
\psdots[dotsize=0.06,dotangle=-270.0](1.705,-0.44996104)
\psdots[dotsize=0.06,dotangle=-270.0](1.905,-0.44996104)
\psdots[dotsize=0.06,dotangle=-270.0](2.105,-0.44996104)
\usefont{T1}{ptm}{m}{n}
\rput(0.29,1.435039){\small $G$}
\psline[linewidth=0.04cm](4.484922,-0.41003904)(4.944922,-0.41003904)
\end{pspicture} 
}

\caption{Construction of $G'$ for a given $G$, notice that not all the edges are shown (for clarity). } %rhochange_2go
\label{fig:setcover}
\end{figure}
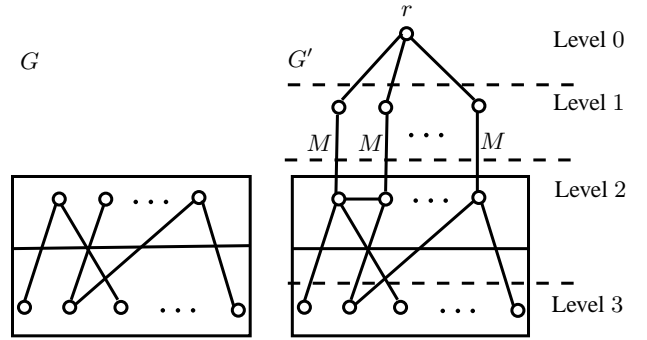

Given an instance, $G$, of the set cover problem, with $k$ sets and $n$ elements, let us construct a new graph $G'$ as follows: Assign a root node $r$, which is the source with the message at the starting time, call this level $0$. Include $k$ nodes in level $1$, representing the $k$ sets in the set cover problem, all connected to the root node, as shown in Figure \ref{fig:setcover} (b). This is followed by the bipartite graph of $G$, which makes up level $2$ and $3$ of $G'$. Connect each of the $k$ nodes in level $2$ to their representative in level $1$ and to all the other nodes in level $2$. Notice the nodes in level $2$ are also connected to their elements in level $3$ of the graph, as shown in the Figure.  We make all the weight on the edges arbitrarily small (say $1$), with the exception of the edges in between the nodes in level $1$ and $2$. We make those edges to be sufficiently large, say $M$, to be specified later.

Assume the the weight on the edges represent the power needed for the message to be transmitted across that edge. If we were to run the optimal DMECB algorithm on $G'$ with $T = 3$ the algorithm would have to act as follows, to be able to cover all the nodes in the given time frame: \newline
\emph{Step $1$:} Root transmits with power $1$, turning on all its $k$ neighbors on level $1$. \newline
\emph{Step $2$:} The algorithm picks a subset of the $k$ nodes on level $1$ to transmit the message. This subset must be chosen to be as small as possible, given the large weight they have to endure to pass on the message on to the bipartite graph, and the fact that DMECB is trying to minimize the total weight. Yet it has to be large enough so that when the nodes in level $2$ transmit, all the nodes in level $3$ would receive the message. The optimal algorithm must be able to find such a subset. \newline
\emph{Step $3$:} The nodes that receive the message in level $2$ transmit the message in this step, turning on all the nodes in level $3$ of the graph, as well as all the nodes in level $2$ of the graph that were not selected for transmission, thus covering the whole graph. \newline
Let us call the solution\footnote{Minimum energy needed for transmission.} of this optimal algorithm $OPT_{DMECB}$. Then the following two lemmas with respect to the above construction of $G$ and $G'$ hold:

\begin{lemma} \label{nplem1}${OPT_{DMECB} \leq M.OPT_{SC} + 1 + OPT_{SC}}$\end{lemma}
\begin{proof} Consider an instance of SC (with graph $G$), whose optimal solution is $OPT_{SC}$. Construct a graph $G'$, as explained and run the DMECB algorithm to get $OPT_{DMECB}$. The above inequality holds by construction of the graph.
\end{proof}

\begin{lemma} \label{nplem2}${OPT_{SC} \leq \frac{OPT_{DMECB}}{M}}$\end{lemma}
\begin{IEEEproof} Consider an instance of DMECB on $G'$ and its optimal solution $OPT_{DMECB}$ for delay $T = 3$. Notice that if $T > 3$, we add additional single nodes (as virtual roots) to reduce the problem to the case where $T=3$. Looking at $G'$, we observe that to meet the delay constraint, by end of step $i$, at least one node in level $i$ must have heard the message - else it is impossible to get the message through to the rest of the levels in the time frame left. Let's say the root is on level $0$. Consider the subset of level $1$ that has come on at the end of time $1$, $\mathbf{s_1}$, and from level $2$ consider the set, $\mathbf{s_2}$, that came on at the end of time step $2$. We now want to show that $\mathbf{s_2}$ is a feasible solution for set cover. To do so, we make the following two claims:
\emph{Claim $1$}: Nodes responsible for turning on $\mathbf{s_2}$ must be a subset of $\mathbf{s_1}$.
\emph{Claim $2$}: $\mathbf{s_2}$  is a feasible solution to set cover.
Claim $1$ holds because only nodes that have received the message by the end of time $1$ can transmit the message at time $2$. Not all of them might transmit though, so $\mathbf{s_2}$ is a subset of corresponding nodes in $\mathbf{s_1}$. Claim $2$ is true because if there exists an element in level $3$ that is not a corresponding node to any node in $\mathbf{s_2}$, it cannot decode by $T=3$.
Therefore, $\mathbf{s_2}$, is a feasible solution to set cover. $OPT_{DMECB}$ must spend at least $M$ for each element of $\mathbf{s_2}$ to come on, so ${OPT_{SC} \leq \frac{OPT_{DMECB}}{M}}$.
\end{IEEEproof}
\begin{theorem} \label{thm.npcomptB}The DMECB problem is $o(\log(n))$ inapproximable, for $T \geq 3$. \end{theorem}

\begin{IEEEproof}  For an instance of the set cover problem, with $k$ being the total number of sets, lemma \ref{nplem1} can be re-written as ${OPT_{DMECB} \leq M.OPT_{SC} + 1 + k}$. We also know by lemma \ref{nplem2} that ${OPT_{SC} \leq \frac{OPT_{DMECB}}{M}}$. Therefore, for a sufficiently large $M$, we can write $OPT_{SC} = \frac{OPT_{DMECB}}{M} + o(1)$. Therefore, the reduction used in construction of the graph $G'$ preserves the approximation factor. That is, if one can find an $\alpha$-approximation for DMECB, by extension there must exist an $\alpha$-approximation for set cover. We know, by \cite{Raz97}, that the set cover problem is $o(\log(n))$ inapproximable, thus DMECB must be $o(\log(n))$ inapproximable. In other words, finding a polynomial time approximation algorithm that approximates $OPT_{DMECB}$ with a factor of $o(\log(n))$ is NP-hard.\end{IEEEproof}

The DMECB problem can be solved in polynomial time for cases when $T< 3$. The optimal algorithm for $T=1$ is trivial and an optimal polynomial algorithm for $T=2$ is discussed in section \ref{sec:optBroadcast}\footnote{DMECT\_{go} algorithm, discussed in section \ref{sec:optBroadcast}, along with an ordering based on channel gains from the source, provides an optimal polynomial time algorithm for DMECB for the case when $T=2$.}.  It is also trivial to verify the feasibility of a given power allocation, and verify whether or not it satisfies the decision version of DMECB given in section \ref{sec:formul}. Therefore, the problem belongs to the class of NP. Notice that the inapproximability result, given by Theorem \ref{thm.npcomptB}, is stronger than, and implies, the NP-completeness result. It is also worth noticing that without any delay constraint (i.e. when $T\geq n$), the problem is still NP-complete and the proof can be obtained, using directed Hamiltonean path, following the approach in \cite{Maric04}.

\subsection{Inapproximablity of DMECM}

The proof of the following theorem, follows from Theorem \ref{thm.npcomptB} by noticing that broadcast can be thought of as a special case of multicast.

\begin{theorem} The DMECM problem is $o(\log(n))$ inapproximable, for $T \geq 3$. \end{theorem}

\subsection{Hardness Results for DMECU}

In the unicast case, the hardness of the problem depends on whether we are using EA or MIA. In the former case, DMECU can be shown to be polynomially solvable and the algorithm for that is provided in section \ref{sec:unicast}. In the remainder of this section, we discuss DMECU with MIA. 

\begin{figure}[!h]
\center
% Generated with LaTeXDraw 2.0.8
% Mon Dec 20 02:35:10 PST 2010
% \usepackage[usenames,dvipsnames]{pstricks}
% \usepackage{epsfig}
% \usepackage{pst-grad} % For gradients
% \usepackage{pst-plot} % For axes
\scalebox{1} % Change this value to rescale the drawing.
{
\begin{pspicture}(0,-1.92)(6.7003126,1.92)
\pscircle[linewidth=0.04,dimen=outer](0.5128125,0.36){0.1}
\pscircle[linewidth=0.04,dimen=outer](1.4928125,1.26){0.1}
\pscircle[linewidth=0.04,dimen=outer](1.4928125,0.64){0.1}
\pscircle[linewidth=0.04,dimen=outer](1.4728125,-0.6){0.1}
\psline[linewidth=0.04cm](0.5128125,0.42)(1.3928125,1.22)
\psline[linewidth=0.04cm](0.5728125,0.38)(1.4128125,0.62)
\psline[linewidth=0.04cm](0.5728125,0.3)(1.3928125,-0.54)
\pscircle[linewidth=0.04,dimen=outer](2.8528125,1.26){0.1}
\pscircle[linewidth=0.04,dimen=outer](2.8528125,0.64){0.1}
\pscircle[linewidth=0.04,dimen=outer](2.8328125,-0.6){0.1}
\psline[linewidth=0.04cm](1.5728126,1.28)(2.7728126,1.3)
\psline[linewidth=0.04cm](1.5728126,0.62)(2.7528124,0.64)
\psline[linewidth=0.04cm](1.5528125,-0.58)(2.7928126,-0.58)
\pscircle[linewidth=0.04,dimen=outer](4.2928123,1.72){0.1}
\pscircle[linewidth=0.04,dimen=outer](4.2928123,1.12){0.1}
\pscircle[linewidth=0.04,dimen=outer](4.3328123,-1.12){0.1}
\pscircle[linewidth=0.04,dimen=outer](4.2928123,0.44){0.1}
\psline[linewidth=0.04cm](2.9128125,1.3)(4.2128124,1.72)
\psline[linewidth=0.04cm](2.9128125,1.24)(4.1928124,0.48)
\psline[linewidth=0.04cm](2.9128125,0.66)(4.1928124,1.1)
\psline[linewidth=0.04cm](2.8728125,-0.64)(4.2728124,-1.06)
\psline[linewidth=0.04cm](2.8728125,-0.54)(4.2528124,1.04)
\pscircle[linewidth=0.04,dimen=outer](5.7928123,0.36){0.1}
\psline[linewidth=0.04cm](4.3528123,1.7)(5.7528124,0.44)
\psline[linewidth=0.04cm](4.3528123,1.1)(5.6928124,0.38)
\psline[linewidth=0.04cm](4.3528123,0.44)(5.7128124,0.34)
\psline[linewidth=0.04cm](4.3928127,-1.1)(5.7528124,0.3)
\psdots[dotsize=0.06](2.0028124,0.3)
\psdots[dotsize=0.06](2.0028124,0.1)
\psdots[dotsize=0.06](2.0028124,-0.1)
\psdots[dotsize=0.06](3.5728126,-0.16)
\psdots[dotsize=0.06](3.5728126,-0.36)
\psdots[dotsize=0.06](3.5728126,-0.56)
\psdots[dotsize=0.06](4.7928123,0.1)
\psdots[dotsize=0.06](4.7928123,-0.1)
\psdots[dotsize=0.06](4.7928123,-0.3)
\usefont{T1}{ptm}{m}{n}
\rput(5.05,1.49){$M$}
\usefont{T1}{ptm}{m}{n}
\rput(3.4996874,1.69){$1$}
\usefont{T1}{ptm}{m}{n}
\rput(2.0189064,1.52){$W$}
\usefont{T1}{ptm}{m}{n}
\rput(0.7796875,0.95){$1$}
\usefont{T1}{ptm}{m}{n}
\rput(0.2128125,0.35){$r$}
\usefont{T1}{ptm}{m}{n}
\rput(6.117344,0.41){$d$}
\psline[linewidth=0.04cm,linestyle=dashed,dash=0.16cm 0.16cm](1.1928124,1.9)(1.1928124,-1.9)
\psline[linewidth=0.04cm,linestyle=dashed,dash=0.16cm 0.16cm](2.5928125,1.9)(2.5928125,-1.9)
\psline[linewidth=0.04cm,linestyle=dashed,dash=0.16cm 0.16cm](3.9928124,1.9)(3.9928124,-1.9)
\psline[linewidth=0.04cm,linestyle=dashed,dash=0.16cm 0.16cm](5.3928127,1.9)(5.3928127,-1.9)
\usefont{T1}{ptm}{m}{n}
\rput(0.5228125,-1.59){Level $0$}
\usefont{T1}{ptm}{m}{n}
\rput(1.9071875,-1.59){Level $1$}
\usefont{T1}{ptm}{m}{n}
\rput(3.3217187,-1.59){Level $2$}
\usefont{T1}{ptm}{m}{n}
\rput(4.7142186,-1.59){Level $3$}
\usefont{T1}{ptm}{m}{n}
\rput(6.1229687,-1.59){Level $4$}
\end{pspicture}
}
\caption{Construction of $G'$ for a given $G$ in DMECU, notice that not all the edges are shown (for clarity). } 
\label{fig:uniMIA}
\end{figure}
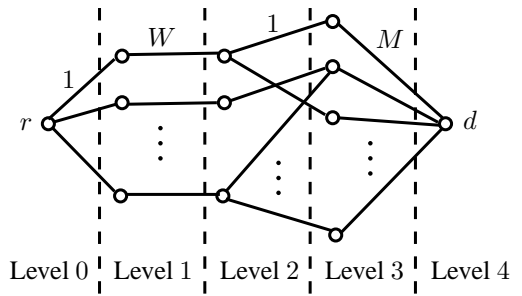

Given an instance, $G$, of the set cover problem, with $k$ sets and $n$ elements, similar to that in section \ref{sec:npcomB}, let us construct a new graph $G'$ as follows: Assign a root node $r$, which is the source with the message at the starting time, call this level $0$. Include $k$ nodes in level $1$, representing the $k$ sets in the set cover problem, all connected to the root node with a small weight (say weight $1$), as shown in Figure \ref{fig:uniMIA}. This is followed by the bipartite graph of $G$, which makes up level $2$ and $3$ of $G'$. Connect each of the $k$ nodes in level $2$ to their representative in level $1$ with edge weights, of say $W$. Notice the nodes in level $2$ are also connected to their elements in level $3$ of the graph, as shown in the Figure, with low-weight edges.  Add a single destination node $d$, in level $4$ and connect all the nodes in level $3$ to $d$. Let the channel between all nodes on level $3$ and destination $d$ be equal and of gain $h$. Therefore, the edge weight on the edges connecting the level $3$ nodes to $d$, can be assigned to be $M$, where $M$ is defined so that the following equality holds: $\log(1+Mh)= \theta$.

Assume that the weight on the edges represent the power needed for the message to be transmitted across that edge. If we were to run the optimal DMECU algorithm on $G'$ with $T = 4$ the algorithm would have to act as follows, to be able to turn node $d$ on within the given time frame: \newline
\emph{Step $1$:} Root transmits with power $1$, turning on all its $k$ neighbors on level $1$. \newline
\emph{Step $2$:} The algorithm picks a subset the $k$ nodes on level $1$ to transmit the message to the nodes in level $2$. \newline
\emph{Step $3$:} A subset of nodes that have received the message in level $2$, transmit the message in this step, turning on a subset of nodes in level $3$ of the graph. \newline
\emph{Step $4$:} A subset of nodes that  have received the message in level $3$, transmit the message in this step with sufficient power to turn on $d$. 
\newline
Let us call the solution of this optimal algorithm $OPT_{DMECU}$.
\begin{theorem} \label{thm.npUMIA}The DMECU problem, with MIA, is NP-complete for $T \geq 4$. \end{theorem}
\begin{IEEEproof}  Given an instance of $G$, we construct $G'$ as above. Let us run DMECU on $G'$ and call the optimal solution $OPT_{DMECU}$ for delay $T = 4$. Notice that if $T > 4$, we add additional single nodes (as virtual roots) to reduce the problem to the case where $T=4$. Define $p$ to satisfy the following:  $n\log(1+ph) = \log(1+Mh)$, meaning $p$ is the power required for nodes on level $3$ to turn on $d$, if all of them were transmitting at the same time. \emph{Claim}:  $OPT_{DMECU}$ needs to use all the nodes in level $3$ for transmission. This claim holds by contradiction, as follows: If all the nodes on level $3$ are used for transmission, each node on that level must transmit with power $p$. Let's assume one of the nodes in that level is not used for transmission. Then the remaining nodes in level $3$ need to transmit with power $p'$, where $n\log(1+ph) = (n-1)\log(1+p'h) = \theta$. Therefore, the ratio of the sum power needed with one fewer node transmitting to the case where all nodes in level $3$ are transmitting can be written as $\frac{(n-1)p'}{np} = \frac{(n-1)(e^{\theta/{n-1}}-1)}{n(e^{\theta/n}-1)}$, for sufficiently large $\theta$, this ratio can become arbitrarily large. Therefore, for sufficiently large $\theta$, the claim holds. Given the claim holds, we know that by definition $OPT_{SC}$ provides the optimal way (minimum energy) to turn on all the nodes in level $3$ within the required time frame, therefore, for non-zero edge weights $OPT_{DMECU}$ needs to optimally solve the set cover problem in step $2$. 
\end{IEEEproof}

It is worth noticing that all the hardness results presented in this section extend to the case where there is no memory.

 %NP, log(n), lowerbounds
\section{OPTIMAL Transmission GIVEN ORDERING\label{sec:optBroadcast}}

In Section \ref{sec:npcomp}, we proved NP-hardness for several variations of the DMECT problem (namely, the DMECB, DMECM  and DMECU (with MIA), with the former two being $o(\log(n))$ innaproximable). In this section, we break this NP-hard problem into three subproblems, namely ordering, scheduling and power allocation, and we propose an optimal polynomial time algorithm for joint scheduling and power allocation when the ordering is given. We evaluate a heuristic for the ordering in Section \ref{sec:simul}.

\begin{definition} An ordering, for a vector of $n$ nodes, is an array of indices from $1$ to $n$; any node that has decoded the message will only be allowed to retransmit when all nodes with smaller index have also decoded the message (and are thus allowed to take part in transmission).
\end{definition}

Given an ordering, what remains to be determined is which nodes should take part in transmission, how much power they should transmit with and at what time slots, such that minimum energy is consumed while delay constraints are satisfied.

\subsection{Instantaneous optimal power allocation}

If we know which nodes are transmitting the message and which nodes are receiving it, at any single time-slot, we can use a convex program (CP) to determine the optimal power allocation for that time slot. Consider an \emph{ordered} vector of $n$ nodes $(1,...,k,...,i,...,n)$. Let us assume that by time slot $t$, node $1$ to $i$ have decoded the message and nodes $i+1$ to $n$ are to decode it during that time slot. At time instance $t$, the optimal instantaneous power allocation for a set of transmitting nodes (say $S(t)=(k,...,i)$) to turn on a set of receiving nodes nodes (say $R(t) =(i+1,...,n)$) can be calculated by the following CP:

%\BE{\min \sum_{s \in S(t)} p_{st}}\EE
%such that
%\BE{ y_{rt} \geq \theta~~~~~~~~~~~~~~~\forall r \in R(t) }  \nonumber\EE
%

\begin{eqnarray} \min & \sum_{s \in S(t)} p_{st} \\ s.t.&~~~p_{st} \geq 0,~~\forall s \nonumber  \\ &\hspace{0.5in}y_{rt} \geq \theta,~~\forall r \in R(t)\nonumber  \label{poweralloc} \end{eqnarray}

We use the notation $CP([\{k...i\}, \{i+1... n\}],\theta, H)$ to refer to solution of the above CP. As a notation, $CP([\{x ... y\}, \{z... \alpha\}], \theta, H) = 0$, if $z \geq \alpha$. Notice that in the case where EA is used, this CP simply reduces to a linear program, using the manipulation highlighted in footnote $2$ in section \ref{sec:formul}.

\subsection{Joint Scheduling and power allocation}

Knowing the instantaneous optimal power allocation given the set of senders and receivers at each time slot, all that remains to be done is to determine these sets at each time slot, in order to minimize the overall power while meeting the delay constraint.

Let $C(j,t)$ be the minimum energy needed to cover up to node $j$ in $t$ steps or less. We can calculate this, using the following algorithm:

\BE \label{eq.twodimpoly}{C(j,t) = \min_{k \in (1,..,j)} \left[C(k,t-1) + CP(\{1 ... k\}, \{k+1... j\},\theta, H)\right]}\EE

where $C(k,1) = CP(1,\{2 ... k\},\theta, H)$, $C(1,t) = 0~~~\forall t$.%, and $C(2,t) = h^{-1}_{12}~~~\forall t$. 

Thus, in DMECB, the total minimum cost for covering $n$ nodes by time $T$ can be found by calculating $C(n,T)$. In DMECM, and DMECU (MIA), the same approach could be used, except for node $n$ being replaced by the highest order destination in the former and by the destination order in the latter.

 A pseudocode for the algorithm is presented below:

\begin{algorithm}[htb]
  \caption{Delay constrained minimum energy cooperative transmission, given an ordering (DMECT\_go)}
  \label{alg:sc}
  \begin{algorithmic}[1]
   \STATE \textbf{INPUT:} an ordered array of nodes of size $n$ (where node $i$ is the $i$th node in the array), $T$ (delay), $d$ (destination), $H$ (channel), $\theta$ (threshold).
    \STATE \textbf{OUTPUT:} $C$ (cost matrix)
    \STATE \textbf{Begin:}

    \FOR{$i := 2$ to $n$}
    \STATE $C(i,1) := CP([1, \{2... i\}], \theta, H)$; \ENDFOR
    \FOR{$t := 1$ to $T$}
    \STATE $C(1,t) := 0;$
    %\STATE $C(2,t) :=h^{-1}_{12}$
    \ENDFOR

    \FOR{$t:=2$ to $T$}

    \FOR{$i:=1$ to $d$}

    \FOR{$k:=1$ to $i$}

    \STATE $x(k):= C(k,t-1) + CP([\{1... k\}, \{k+1... i\}], \theta, H)$;

    \ENDFOR
    \STATE $C(i,t) := \min \textbf{x}$

    \ENDFOR
    \ENDFOR

  \end{algorithmic}
\end{algorithm}

The complexity of the optimal scheduling and power allocation can be obtained by inspection of the above algorithm: it invokes at most $O({n^2}T)$ calls to the CP solver, each of which takes polynomial time. Hence the DMECT\_go algorithm that does joint scheduling and power-control is a polynomial time algorithm.

Note that the delay constraint $T$ can be made sufficiently large ($\geq n)$, or removed entirely from the formulation, to cover the case of no delay constraints. In that case the two-dimensional dynamic program proposed in (\ref{eq.twodimpoly}), reduces to a one-dimensional dynamic program:
\BE\label{eq.onedimpoly}{C(n) = \min_k\left[C(k) + CP(\{1 ...k\},\{ k+1 ... n\}),\theta, H\right]}\EE
where $C(n)$ is the minimum cost of covering node $n$ using our cooperative memoryless approach, starting from node $1$ and $C(1) =0$.

  %time constraint, no time-constraint
\section{OPTIMAL UNICAST WITH ENERGY ACCUMULATION} \label{sec:unicast}

In this section we propose an optimal polynomial time algorithms for solving the unicast problem with EA.

\begin{theorem} \label{thm.dmecu} In DMECU with EA, there exists a solution consisting of a simple path between source and destination, which is optimum.
\end{theorem}
\begin{proof} Let us prove by induction:
 For delay $T=1$, the claim is trivially true, as the optimal solution is a direct transmission from the source, $s$, to the given destination, $d$. For $T>1$, we prove the claim by induction. Assume that the claim is true for $T=k-1$. Pick any node in the network as the desired destination $d$. If the message can be transmitted from source $s$ to $d$ with minimum energy in a time frame less than $k$, then an optimal simple path exists by the induction assumption. So consider the case when it takes exactly $T=k$ steps to turn on $d$. The system is memoryless, so $d$ must decode by accumulating the energy transmitted from a set of nodes, $\textbf{v}$, at time $k$. This can be represented as $\log(1+\sum_{v_i\in \textbf{v}} p_{{v_i}k}h_{dv_i} )\geq \theta$. We observe that there must exist a node $v_o \in \textbf{v}$ whose channel to $d$ is equal or better than all the other nodes in $\textbf{v}$. Therefore, given $h_{dv_o} \geq h_{dv_i}, \forall v_i \in \textbf{v}-\{v_o\}$ then $\log(1+\sum_{v_i\in \textbf{v}} p_{{v_i}k}h_{dv_o}) \geq \log(1+\sum_{v_i\in \textbf{v}} p_{{v_i}k}h_{dv_i}) \geq \theta$. In other words, if we add the power from all nodes in $\textbf{v}$ and transmit instead from $v_o$, our solution cannot be worse. $v_o$ must have received the message by time $k-1$, to be able to transmit the message to $d$ at time $k$. We know by the induction assumption that the optimal simple path solution exists from source to any node to deliver the message within $k-1$ time frame. Thus, for $T=k$, there exists a simple path solution between $s$ and $d$, which is optimum.
\end{proof}
Notice that the above theorem holds in the case where there is no delay constraint as well. The proof follows an straightforward modification of the above proof and is omitted for brevity.

\begin{corollary} \label{cor:uni} The Dijkstra's shortest path algorithm provides the optimal ordering in the case of minimum energy memoryless cooperative unicast, when there is no delay constraint.
\end{corollary}
\begin{proof} We have already established that an optimal minimum energy solution exists between source and destination, which is a simple path. The well-known Dijkstra's  shortest path algorithm can find the minimum cost simple path between source and destination. Therefore, Dijkstra's algorithm provides the optimal ordering.
\end{proof}
Using theorem \ref{thm.dmecu} we know that the optimal unicast solution from source to any destination in DMECU (with EA) is given by a simple path. The cost paid by the optimal solution can be calculated using the following algorithm:
Let  $C(i,t)$ be the minimum cost it takes for source node $s$ to turn on $i$, possibly using relays, within at most $t$ time slots. Then we can write: 
\BE{ C(i,t) =\min_{k \in Nr(i)}\left[C(k,t-1)+w(k\rightarrow i)\right] }\EE with $C(s,t) = 0$, for all $t$ and  $C(i,1) = w(s\rightarrow i)$, where  $Nr(i)$ is the set that contains $i$ and its neighboring nodes that have a non-zero channel to $i$, $w(k\rightarrow i)$ represents the power it takes for $k$ to turn on $i$ using direct transmission. Given that, the solution to $OPT_{DMECU}$ (with EA) is given by $C(d,T)$. Computing this lower-bound incurs a running time of $O(n^3)$. 

The unicast case (with EA), with no delay constraint, is still polynomially solvable. Given Theorem \ref{thm.dmecu}  and Corollary \ref{cor:uni}, the optimal solution is simply the weight of the shortest path given by the Dijkstra's algorithm.

It is worth noticing that the crux of the difference between DMECU with EA and with MIA, that allows the former to be polynomially solvable, while the latter is NP-complete, lies in the optimality of single-node transmission. Namely, in the EA case, the multi-transmitters single-receiver case (multi-single) makes no sense as explained above and instead it is optimal to put the combined power into the best channel. This allows for the overall solution to be a simple path. However, in the MIA case, the many-to-one transmission case does in fact make sense. That is due to the property of the $\log$ function, creating an effect similar to what we observe in \emph{water-filling}, where it is best to transmit from the best channel up until some point, then from the second best channel and so forth.

\section{APPROXIMATION ALGORITHM FOR BROADCAST WITH ENERGY ACCUMULATION\label{sec:approx}}

In Section \ref{sec:npcomp}, we proved that DMECB is NP-complete and $o(\log(n))$ inapproximable, therefore it is \emph{hard} to approximate DMECB to a factor strictly better
than $\log(n)$. It is of theoretical interest to know how close we can get to the optimal solution, using a polynomial-time algorithm. In this section we show that existing approximation algorithms for the \emph{bounded-diameter directed Steiner tree problem} can be used to provide an $O(n^{\epsilon})$ approximation for DMECB in the case where EA is used. We do so by proposing an approximation-preserving reduction to the directed Steiner tree problem.

The Steiner tree problem is a classic problem in combinatorial optimization \cite{Vazi01}. We focus on a variation of this problem, namely bounded diameter directed Steiner tree, defined as follows. Given a directed weighted graph $G(V,E)$, a specified root $r \in V$, and a set of \emph{terminal nodes} $X \subseteq V$ ($|X| = n$), the objective is to find the minimum cost arborescence rooted in $r$ and spanning all vertices in $X$, subject to a maximum diameter $T$. Diameter refers to the maximum number of edges on any path in the tree. Notice that the tree may include vertices not in $X$ as well, these are known as \emph{Steiner nodes}. Directed Steiner tree problem is known to be $NP$-complete and $O(\log(n))$ inapproximable \cite{Vazi01}. In \cite{Moses99}, the authors give the first non-trivial approximation algorithms for Steiner tree problems and propose approximation algorithms that can achieve an approximation factor of $O(n^{\epsilon})$  for any fixed $\epsilon > 0$ in polynomial time. To the best of our knowledge this is currently the tightest approximation algorithm known for this problem.

In order to reduce a given instance of the DMECB to an instance of the Steiner tree problem, we first restrict DMECB by not allowing many-transmitter-to-one-receiver (many-to-many) transmissions. Notice that in the proof of theorem \ref{thm.dmecu}, we had established that many-to-one transmissions can be replaced with one-to-one transmissions without loss of optimality. Therefore, by not allowing many-to-many transmissions, we are left with one-to-one and one-to-many transmissions. We call this an integral version of DMECB, DMECB-int. The integrality gap of the weighted set cover problem is shown to be $\log(n)$ \cite{Vazi01}; it is straightforward to extend that result to show that DMECB-int also loses a factor of $\log(n)$, compared to optimal DMECB.

\begin{figure}[!h]
\centering
\includegraphics[width=0.3\textwidth] {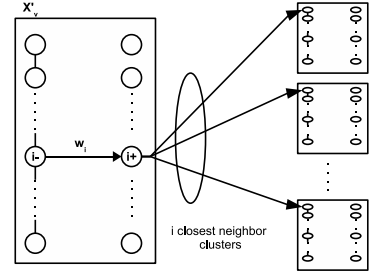}
\caption{A simplified example of how clusters are constructed in $G'$.}
\label{fig:appalg}
 \end{figure}

Consider an instance of DMECB-int, $G(V,E)$, with ($|V| = n$) and $s \in V$ being the source node. To reduce this problem to an instance of directed Steiner tree problem, let us construct a new graph $G'$, consisting of $n$ clusters, $x'$, each corresponding to each node in $G$. Let each cluster be a bipartite graph, with $n$ nodes on the left (marked as $``-"$) and $n$ nodes on the right (marked as $``+"$), as shown in Figure \ref{fig:appalg}. The $``-"$ nodes are intra-connected within a cluster with edges of weight $0$. In each cluster, $x' \in G'$ corresponding to node $v \in G$, the $``+"$ and $``-"$ nodes on each level, $i$, of the bipartite graph are connected to each other with an edge of weight $w_i$, representing the power needed by the corresponding node $v \in G$ to turn on its $i$ closest neighbors. The $i_+$ node is then connected, with edges of weight $0$, to all the $``-"$ nodes in the corresponding neighbor clusters. We further add a single root node, $r \in G'$, and connect it via a zero-weight edge to all the $``-"$ nodes in the cluster corresponding to $s$, $x'_s$. We assign the root $r$ and one desired $``-"$ node from each cluster as terminal nodes and all other nodes in $G'$ as Steiner nodes.

Let us look at an example of this construction, say node $v_1 \in G$, whose closest $3$ neighbors are ($v_2,v_4,v_6$). We have an equivalent cluster $x'_1 \in G'$ corresponding to node $v_1$. $x'_1$ has $2n$ nodes, arranged in $n$ levels. The weight between the two nodes in say level $3$ is equivalent to the power it takes for $v_1$ to turn on ($v_2,v_4,v_6$).  Furthermore, the node $3_+$ in cluster $x'_1$ is connected to the $``-"$ nodes in clusters ($x'_2,x'_4,x'_6$) with edges of weight $0$. This construction allows us to find a way to allow $v_1$ to transmit with different power levels, without knowing what those powers might be in advance. We first add a single root node, $r$, and connect it via a zero-weight edge to all the $``-"$ nodes in the cluster corresponding to $s$.

Run the directed Steiner tree algorithm on $G'$ to obtain a solution. The solution must choose at least one node from each cluster, to meet the mandatory terminal nodes requirement. Recall that each cluster in $G'$ corresponds to a node in $G$ and that multi-multi was not allowed. To convert
the solution of the Steiner tree algorithm on $G'$ to a solution of DMECB-int on $G$, we look at the parent of each cluster, which is a $``+"$ node in another cluster. Let's say we want to see which node turns on $v_6$ by looking at $G'$. We look at the parent of $x'_6$ and see that it's $3_+ \in x'_1$. So in $G$, we figure out that $v_1$ must transmit with enough power to turn on $3$ of its closest neighbor ($w_3$), and it is as a result of this transmission that $v_6$ comes on. Going through all the clusters and their parents, we can establish an ordering and transmission power for all the nodes that should take part in cooperation in $G$, and thus we have a solution for DMECB-int.

\begin{theorem} \label{thm.approx1}For DMECB problem with EA, an $O(n^\epsilon)$ approximation ratio can be achieved in polynomial time, for any fixed $\epsilon > 0$. \end{theorem}
\begin{proof}
As mentioned, the directed Steiner tree is $o(\log(n))$ inapproximable, and the best approximation algorithm currently available \cite{Moses99} gives an $O(n^\epsilon)$ approximation on the optimal solution. We had already lost $O(\log(n))$ to convert DMECB to its integral form. The approximation algorithm proposed in \cite{Moses99} can approximate the optimal integral solution within $O(n^\epsilon)$. Therefore, using the above reduction, and applying the directed Steiner tree approximation algorithm, we can approximate the optimal solution to DMECB within $O(n^\epsilon \times \log(n))$, which is equivalent to $O(n^\epsilon)$.
\end{proof}

The running time of the Steiner tree approximation algorithm is a function of $\epsilon$, and the tighter the approximation, the worse the running time. Similarly, using the above mentioned reduction, the following result holds by directly applying the approximation algorithms in \cite{Moses99}. Detailed discussions of the algorithms in \cite{Moses99} are beyond the scope of this paper. 

\begin{theorem} \label{thm.approx2}For any fixed $T>0$, there is an algorithm which runs in time  $n^{O(T)}$ and gives an  $O(T\log^2(n))$ approximation of the DMECB with EA. \end{theorem}

  %approx_results...steiner tree (log(n) loss + the other part)
\section{PERFORMANCE EVALUATION\label{sec:simul}}

For the simulations, we focus on the broadcast case. We consider a network of $n$ nodes uniformly distributed on a $15$ by $15$ square surface. The transmission starts from a node, arbitrarily located at the left center corner of the network $(0,7)$. The channels between all nodes are static, with independent and exponentially distributed channel gains (corresponding to Rayleigh fading), where $h_ij$ denotes the channel gain between node $i$ and $j$. The mean value of the channel between two nodes, $\overline{h_{ij}}$, is chosen to decay with the distance between the nodes, so that $\overline{h_{ij}} = d_{ij}^{-\eta}$, with $d_{ij}$ being the distance between nodes $i$ and $j$ and $\eta$ being the path loss exponent. The corresponding distribution for the channel gains is then given by
\[f_{h_{ij}}(h_{ij})= {\frac{1}{\overline{h_{ij}}}}\exp \left(\frac{h_{ij}(k)}{\overline{h_{ij}}}\right)\] Based on the intuition developed in section \ref{sec:unicast}, we use the Dijkstra's shortest path algorithm as our ordering heuristic. Simulations are repeated multiple times with the same node locations but different fading realizations and average values are shown in the graphs. Notice that the minimum power calculated by different algorithms, shown on the y-axes of the graphs in this section, are normalized by unit power (rendering it unit-less). The value of $\theta$ is, arbitrarily, chosen to be $\log(2)$ throughout this section.

\begin{figure}[!htb]
\hspace{-0.3cm}
\includegraphics[height=5cm,width=10cm]{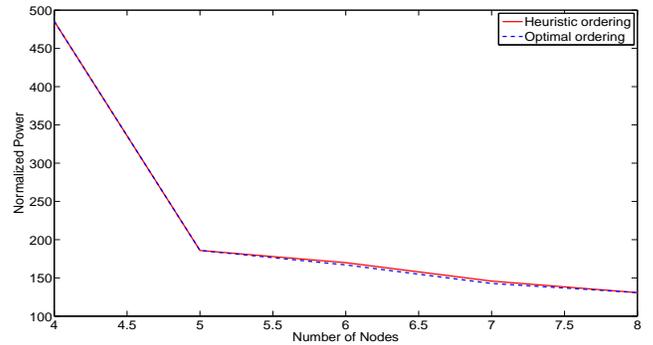} %{Graph1.eps}
\caption{Performance with optimal ordering vs Dijkstra's algorithm-based heuristic ordering.} %rhochange_2go
\label{fig:optiFigure}
\end{figure}

In Figure \ref{fig:optiFigure}, we calculate the optimal ordering by brute-force for a small number of nodes and compare the performance of our algorithm, which uses Dijkstra's shortest path-based ordering, with the optimal performance. The results for the broadcast case, with EA, is shown in this figure. As can be seen, Dijkstra's algorithm provides a good heuristic for ordering in this example and will be used throughout this section. Note that, as can be seen in the figure, although the problem was proved to be $o(\log(n))$ inapproximable, it is possible to achieve near-optimal results in polynomial time in certain practical settings where the network does not have any pathological properties. The inapproximablity results contain all possible (including pathological networks) scenarios.

We next compare the performance of our cooperative algorithm with a smartly designed non-cooperative algorithm (both using EA). Notice that in our cooperative algorithm we make use of the \emph{wireless broadcast advantage (WBA)}, where transmission by one node can be received by multiple nodes and \emph{cooperative advantage}, where a node can accumulate power from multiple transmitters. If an algorithm is using the WBA, but not the cooperative advantage, it can be thought of as an integral version of DMECB. This means that each node can receive the message from one transmitter only (and cannot accumulate from multiple transmitters), however one transmitter can transmit to multiple receivers. We had established in Section VII that DMECB-int is also NP complete. It is however interesting to note that DMECB-int needs to solve a weighted set cover problem when allocating powers as well; we know that set cover problem is $o(\log(n))$ inapproximable \cite{Vazi01}, so the non-cooperative case is $o(\log(n))$ inapproximable, even when ordering is provided. Greedy algorithms exist \cite{Vazi01} that give $O(\log(n))$ approximations for the weighted set cover problem, and thus provide a tight polynomial time approximation. Therefore, to simulate a smart non-cooperative algorithm, we use Dijkstra's algorithm-based ordering and the DMECT algorithm of section \ref{sec:optBroadcast}, with the exception that instead of using an LP we use the greedy algorithm for power allocation.

\begin{figure}[!htb]
\hspace{-0.5cm}
\includegraphics[height=5cm,width=10cm]  {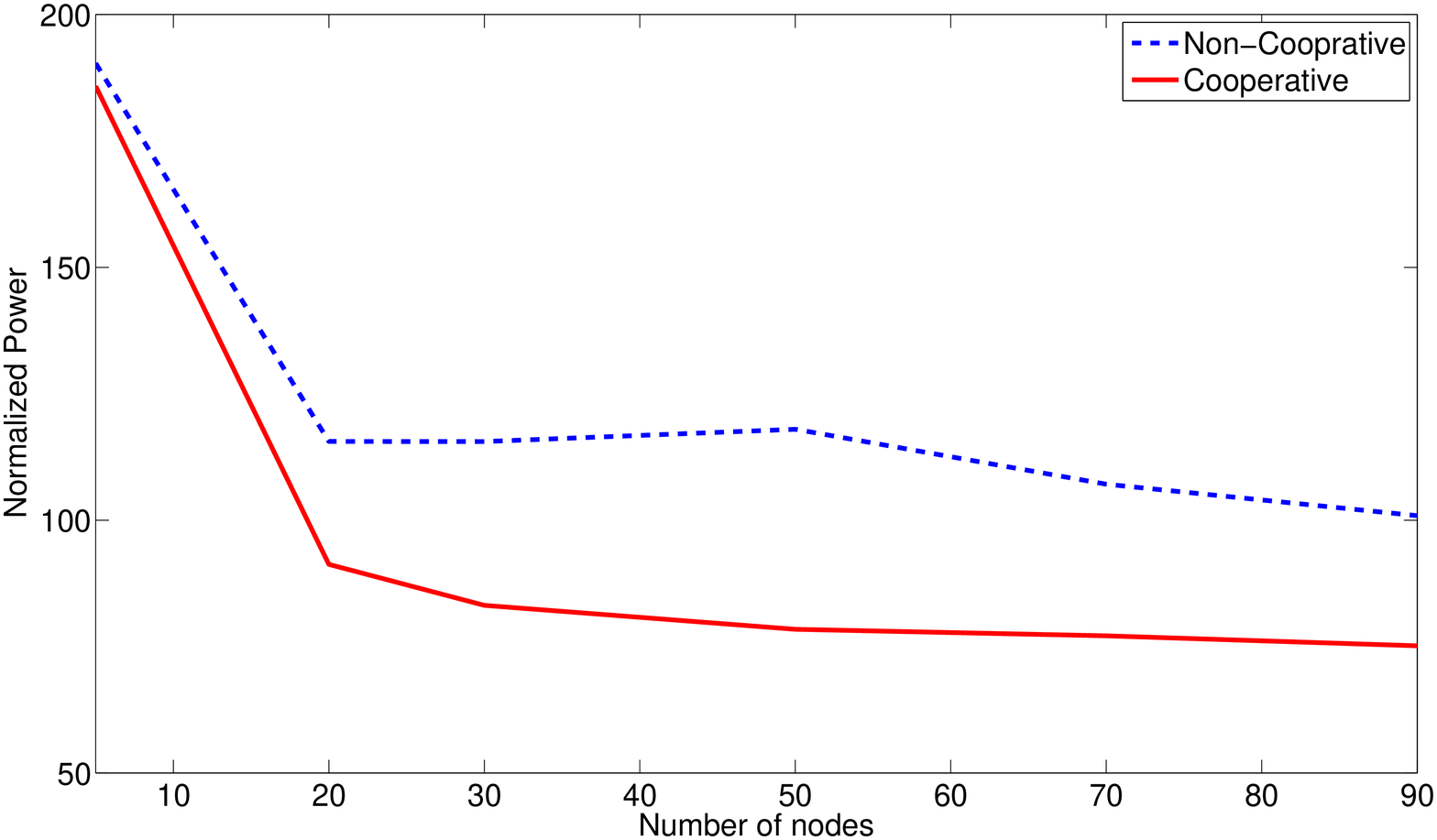}
\caption{Effect of cooperation in broadcast.} %rhochange_2go
\label{fig:coopeffect}
\end{figure}

\begin{figure}[!htb]
\hspace{-0.5cm}
\includegraphics[height=5cm,width=10cm]  {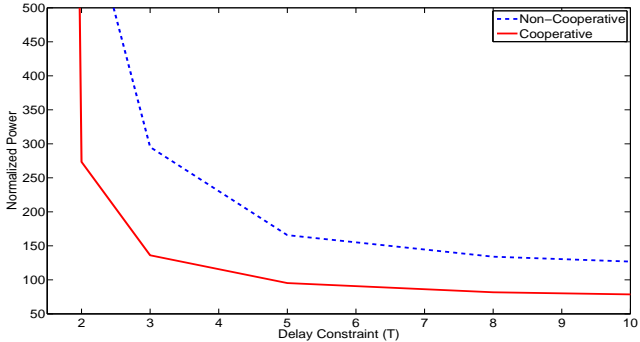}
\caption{Power-delay tradeoff in cooperative vs non-cooperative case.} %rhochange_2go
\label{fig:G3}
\end{figure}

The performance comparison between our proposed cooperative algorithm and the smart non-cooperative algorithm, for different values of $n$ is shown in Figure \ref{fig:coopeffect} and the power-delay tradeoff for cooperative and non-cooperative algorithms are presented in \ref{fig:G3}. As can be seen, the cooperative algorithm outperforms the non-cooperative algorithm, and the advantage is more pronounced when a delay constraint is imposed.

\begin{figure}[!htb]
\hspace{-1cm}
\includegraphics[height=5cm,width=10cm]{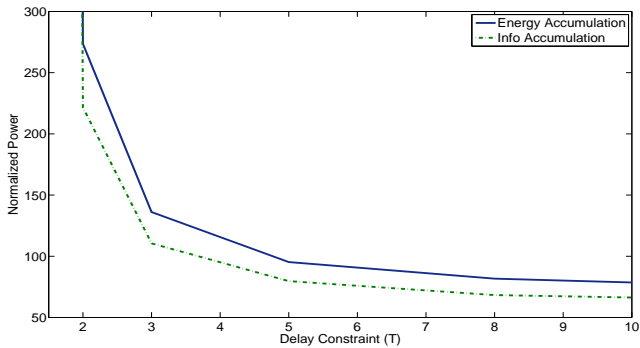}%{info.eps} %{Graph1.eps}
\caption{Energy accumulation vs mutual information accumulation.} %rhochange_2go
\label{fig:EAvsMIA}
\end{figure}

The performance gains obtained by using MIA is shown in Figure \ref{fig:EAvsMIA}, for a sample network of $30$ nodes.

\begin{figure}[!htb]
\hspace{-0.3cm}
\includegraphics[height=5cm,width=10cm] {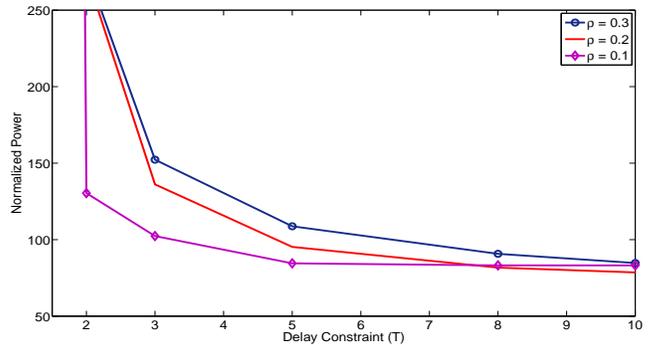}
\caption{Effect of network density on power-delay tradeoff.} %rhochange_2go
\label{fig:G4}
\end{figure}

\begin{figure}[!htb]
\hspace{-0.3cm}
\includegraphics[height=5cm,width=10cm]{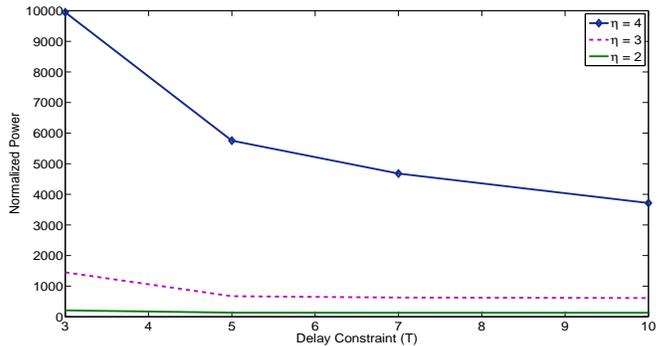}%{Eta_infoAcc_30nodes.eps} %{Graph1.eps}
\caption{Different $\eta$ values for information accumulation.} %rhochange_2go
\label{fig:MIAeta}
\end{figure}

We next study the power-delay tradeoff of the cooperative algorithm for different channel conditions and different values of network density $\rho$ (in nodes/area). Figure \ref{fig:G4} and Figure \ref{fig:MIAeta}, show results for EA and MIA, respectively. These figures highlight the sensitivity of the dense networks and those with poor channel conditions to delay constraints and the importance of having smart algorithms to minimize the energy consumption.

 %\newpage
%%%\input{conclusion}
\section{CONCLUSIONS\label{sec:conc}}

We have formulated in this work the novel problem of delay constrained minimum energy cooperative transmission (DMECT) in memoryless wireless networks, encompassing both EA and MIA. We have shown that this problem is 
$o(\log n)$ inapproximable in broadcast and multicast cases and is NP-complete in the unicast case when mutual information accumulation is used. For the broadcast case with EA, we have presented an approximation algorithm. We also developed a polynomial-time algorithm that can solve
this NP-hard problem optimally for a fixed transmission ordering. Our empirical results suggest that for practical settings, a near-optimal ordering can be obtained by using Dijkstra's shortest path algorithm. We have further showed that the unicast case can be solved optimally and in polynomial time when EA is used. We have studied the energy-delay tradeoffs and the performance gain of MIA using simulations, and evaluated the performance of our algorithm under varying conditions.

The summary of the algorithmic results developed in this paper are presented in Tables \ref{tb:neg} and \ref{tb:pos}. Besides tackling the unsolved problems signified by empty slots in this table, there are a number of other interesting directions for future work. In this paper we have focused on the static problem with full channel-state information, which allows for centralized decision making. Distributed solutions to this problem which would be particularly suitable for more dynamic settings.
From an analytical perspective, deriving inapproximability results for the NP-complete cases in the tables, as well as tighter approximation results are interesting avenues for further research. Evaluating the proposed algorithms under more realistic settings (through
more detailed simulations of physical layer implementation or through direct implementation on software radio platforms, and the use of more realistic energy models) would certainly help in moving this work towards practice. Finally, our work, like most work in this domain of cooperative broadcasts, has focused on the single flow setting. It is of interest to study generalizations that allow for multiple simultaneous flows in the network.

% that's all folks
\end{document}